\documentclass[11pt]{amsart}
\usepackage[margin=1.1in]{geometry}
\usepackage{amsmath,amssymb,amsthm,mathtools}
\usepackage{enumitem}
\usepackage{booktabs}
\usepackage{tikz}
\usepackage{pgfplots}
\pgfplotsset{compat=1.17}
\usepackage[colorlinks=true,linkcolor=blue!50!black,citecolor=green!40!black,urlcolor=blue!50!black]{hyperref}

\newtheorem{theorem}{Theorem}[section]
\newtheorem{proposition}[theorem]{Proposition}
\newtheorem{lemma}[theorem]{Lemma}
\newtheorem{corollary}[theorem]{Corollary}
\newtheorem{maintheorem}{Theorem}

\theoremstyle{definition}

\newtheorem{example}[theorem]{Example}
\theoremstyle{remark}
\newtheorem{remark}[theorem]{Remark}

\newcommand{\N}{\mathbb{N}}
\newcommand{\Z}{\mathbb{Z}}
\newcommand{\R}{\mathbb{R}}
\newcommand{\Q}{\mathbb{Q}}

\newcommand{\Prob}{\mathbb{P}}
\newcommand{\E}{\mathbb{E}}
\newcommand{\W}{W_\beta}
\newcommand{\eps}{\varepsilon}

\makeatletter
\@namedef{subjclassname@2020}{\textup{2020} Mathematics Subject Classification}
\makeatother

\title{The line crossing problem for biased monotonic random walks}
\author{Mohammad Javaheri}
\address{Department of Mathematics and Statistics, Siena University\\
515 Loudon Road, Loudonville, NY 12211, USA}
\email{mjavaheri@siena.edu}
\date{\today}
\subjclass[2020]{60G50, 05A15, 60J80, 05A10}
\keywords{Biased random walk, lattice paths, Fuss--Catalan numbers, ballot problem, kernel method, Galton--Watson process, ruin probability}

\begin{document}

\begin{abstract}
A $\beta$-biased monotonic random walk in the plane moves from $(x,y)$ to $(x+1,y)$ or $(x,y+1)$ with probabilities $1/(\beta+1)$ and $\beta/(\beta+1)$. We study $\Phi_\beta(\alpha,d)$, the probability that the walk started at the origin ever visits a lattice point strictly above the line $y=\alpha x+d$, and $\Psi_\beta(\alpha,d)$, the probability that it visits a lattice point \emph{on} the line after leaving the origin. For integer slopes $\alpha=p$ we show $\Phi_\beta(p,d)=\rho^{\lfloor d\rfloor+1}$, where $\rho$ is the smallest positive root of $y^{p+1}-(\beta+1)y+\beta=0$, give the Fuss--Catalan series for $\rho$, identify $\rho$ with the extinction probability of a Galton--Watson process, and compute the exact law and mean of the first crossing point. We prove that the return probability is $\Psi_\beta(p,0)=(p+1)(1-\beta/((\beta+1)\rho))$, which equals $(p+1)/(\beta+1)$ for $\beta\ge p$. For rational slopes $\alpha=r/s$ we obtain $\Phi_\beta$ in closed form from the $s$ large roots of the kernel $\tfrac{\beta}{\beta+1}z^{r+s}-z^r+\tfrac1{\beta+1}$, and $\Psi_\beta$ as a ratio of Green's functions. For arbitrary real $\alpha$ we prove that the walk crosses every line $y=\alpha x+d$ almost surely if and only if $\beta\ge\alpha$, and we establish the two-sided bound $e^{-\theta(d+1)}\le\Phi_\beta(\alpha,d)<e^{-\theta d}$ with an explicit Lundberg exponent $\theta$. Finally, $\alpha\mapsto\Phi_\beta(\alpha,0)$ is a pure jump function: it is constant on $[0,\beta]$, continuous at every irrational point and has a jump at every rational $\alpha>\beta$.
\end{abstract}

\maketitle

\section{Introduction}\label{sec:intro}

\subsection{Motivation}
It is classical that the simple symmetric random walk on $\Z$ visits every integer infinitely often almost surely; equivalently, P\'olya's return probability in dimension one equals $1$ \cite{Polya,Finch}. Chung and Fuchs \cite{ChungFuchs} extended this (see also \cite{Domb}) to sums of i.i.d.\ random variables, and the recurrence/transience dichotomy has since become a cornerstone of probability \cite{Spitzer,Durrett}.

A one-dimensional $\pm1$ walk is the same thing as a \emph{monotonic} lattice path in the plane: record an up-step $(0,1)$ for every $+1$ and a right-step $(1,0)$ for every $-1$. The level of the one-dimensional walk is then $y-x$, and returning to $0$ means touching the diagonal $y=x$. This two-dimensional picture suggests a natural family of questions in which the diagonal is replaced by an arbitrary line $y=\alpha x+d$ and the fair coin by a biased one. We say that a random walk in the plane is \emph{$\beta$-biased monotonic} ($\beta>0$) if from $(x,y)$ it moves to
\[
(x+1,y)\ \text{ with probability } a:=\frac{1}{\beta+1},\qquad
(x,y+1)\ \text{ with probability } b:=\frac{\beta}{\beta+1},
\]
independently of the past. We write $\W=(X_k,Y_k)_{k\ge0}$ for this walk started at the origin. The parameter $\beta=b/a$ is the expected number of up-steps per right-step, so the walk travels asymptotically along the ray $y=\beta x$.

The term ``biased random walk'' is borrowed from cell biology, where the run-and-tumble motion of bacteria such as \emph{E.~coli} is modelled by a random walk whose step distribution is tilted by a chemical gradient \cite{BergBrown,Codling}. The monotonic model is the simplest caricature of such motion: in each unit of time the organism either advances (an up-step) or is diverted (a right-step), and one asks whether it ever reaches a moving target whose position grows linearly in time.

The questions we study are:
\begin{itemize}[leftmargin=2em]
\item[$\Phi$:] What is the probability $\Phi_\beta(\alpha,d)$ that $\W$ \emph{crosses} the line $y=\alpha x+d$, i.e.\ visits some lattice point $(x,y)$ with $y>\alpha x+d$?
\item[$\Psi$:] What is the probability $\Psi_\beta(\alpha,d)$ that $\W$ \emph{meets} the line, i.e.\ visits a lattice point on it, after leaving the origin?
\end{itemize}

Besides their intrinsic interest, these probabilities appear in several guises.
\begin{enumerate}[leftmargin=2em]
\item \emph{Ballot problems.} When $\alpha=p\in\N$, the paths that stay weakly below $y=px$ are counted by the Fuss--Catalan numbers \cite{Fuss,HiltonPedersen,Stanley}, and $\Phi_\beta$ is a probabilistic weighting of the generalized ballot problem \cite{Takacs,Bizley,BanderierFlajolet}.
\item \emph{Branching processes.} Via the \L ukasiewicz correspondence, $\Phi_\beta(p,d)$ is the probability that a Galton--Watson process with offspring law $\Prob(0)=b$, $\Prob(p+1)=a$, started from $d+1$ individuals, becomes extinct (third proof of Theorem~\ref{thm:A} in \S\ref{sec:integer}).
\item \emph{Risk theory.} $\Phi_\beta(p,d)$ is the ruin probability in the discrete-time \emph{dual} risk model with initial capital $d$, unit cost per period and gains of size $p+1$ arriving with probability $a$ per period \cite{Gerber88,AGS07} (Remark~\ref{rem:risk}); for general $\alpha$ our two-sided estimate (Theorem~\ref{thm:D}) is a sharp form of the Cram\'er--Lundberg inequality \cite{Asmussen}.
\end{enumerate}

\subsection{Main results}
Let $p\ge0$ be an integer. Denote by $\rho=\rho_{\beta,p}$ the smallest positive root of
\begin{equation}\label{eq:rho}
y^{p+1}-(\beta+1)y+\beta=0,
\end{equation}
and let $M(p,n)=\frac{1}{pn+1}\binom{(p+1)n}{n}$ be the Fuss--Catalan numbers.

\begin{maintheorem}[Crossing a line of integer slope]\label{thm:A}
Let $\beta>0$, $p\in\{0,1,2,\dots\}$ and $d\ge0$ real. Then
\[
\Phi_\beta(p,d)=\rho^{\lfloor d\rfloor+1},\qquad
\rho=\sum_{n\ge0}M(p,n)\,a^n b^{pn+1}.
\]
Moreover $\rho=1$ if and only if $\beta\ge p$. The first crossing occurs at a point $(n,pn+\lfloor d\rfloor+1)$ with
\[
\Prob\bigl(\text{first crossing at column } n\bigr)=\frac{k}{(p+1)n+k}\binom{(p+1)n+k}{n}a^n b^{pn+k},\qquad k=\lfloor d\rfloor+1,
\]
and if $\beta>p$, the expected column of the first crossing is $k/(\beta-p)$.
\end{maintheorem}

\begin{maintheorem}[Meeting a line of integer slope]\label{thm:B}
Let $\beta>0$ and $p,d\in\{0,1,2,\dots\}$. Then $\Psi_\beta(p,d)=\rho^d$ for $d\ge1$, and
\[
\Psi_\beta(p,0)=(p+1)\Bigl(1-\frac{\beta}{(\beta+1)\rho}\Bigr)=(p+1)\,z_*,
\]
where $z_*$ is the smallest positive root of $z(1-z)^p=a b^p$. In particular $\Psi_\beta(p,0)=(p+1)/(\beta+1)$ for $\beta\ge p$, and $\Psi_\beta(p,0)=1$ if and only if $\beta=p$. Equivalently, the number of monotonic lattice paths from $(0,0)$ to $(n,pn)$, $n\ge1$, that meet the line $y=px$ only at their endpoints is
\[
F(p,n)=\frac{p+1}{n}\binom{(p+1)n-2}{n-1}.
\]
\end{maintheorem}

\begin{maintheorem}[Rational slopes]\label{thm:C}
Let $\alpha=r/s$ with $r,s\ge1$ coprime and $\beta<\alpha$. The kernel $K(z)=bz^{r+s}-z^r+a$ has exactly $s$ roots $v_1,\dots,v_s$ in $|z|>1$, and
\[
\Phi_\beta(r/s,d)=1-[z^{\lfloor sd\rfloor}]\;\frac{1}{1-z}\prod_{j=1}^s\frac{v_j-1}{v_j-z},
\qquad\text{in particular}\qquad
\Phi_\beta(r/s,0)=1-\prod_{j=1}^s\Bigl(1-\frac1{v_j}\Bigr).
\]
For instance $\Phi_\beta(1/s,0)=\min\{1,s\beta\}$. For every $\beta\ne r/s$,
\[
\Psi_\beta(r/s,0)=1-\Bigl(\sum_{k\ge0}\binom{(r+s)k}{sk}\Bigl(\frac{\beta^r}{(\beta+1)^{r+s}}\Bigr)^k\Bigr)^{-1}.
\]
\end{maintheorem}

\begin{maintheorem}[Arbitrary slopes]\label{thm:D}
Let $\alpha\ge0$ and $d\ge0$ be real.
\begin{enumerate}[label=\textup{(\roman*)}]
\item $\Phi_\beta(\alpha,d)=1$ for all $d\ge0$ if and only if $\beta\ge\alpha$; if $\beta<\alpha$ then $\Phi_\beta(\alpha,d)<1$ for every $d\ge0$.
\item If $\beta<\alpha$, let $\theta>0$ be the unique positive root of $b e^{\theta}+a e^{-\alpha\theta}=1$. Then
\[
e^{-\theta(d+1)}\le\Phi_\beta(\alpha,d)<e^{-\theta d},
\]
with equality on the left when $\alpha$ and $d$ are integers.
\item The function $\alpha\mapsto\Phi_\beta(\alpha,0)$ equals $1$ on $[0,\beta]$, is right-continuous and strictly decreasing on $[\beta,\infty)$, tends to $b$ as $\alpha\to\infty$, is continuous at every irrational $\alpha$, and has a jump at every rational $\alpha>\beta$. It is a pure jump function on $(\beta,\infty)$.
\end{enumerate}
\end{maintheorem}

Figure~\ref{fig:phi} illustrates Theorem~\ref{thm:D}.

\begin{figure}[t]
\centering
\begin{tikzpicture}
\begin{axis}[
  width=0.92\textwidth, height=7.2cm,
  xmin=0, xmax=5, ymin=0.45, ymax=1.03,
  xlabel={slope $\alpha$}, ylabel={$\Phi_\beta(\alpha,0)$},
  legend cell align=left,
  legend style={font=\small, at={(0.5,-0.22)}, anchor=north, legend columns=1, draw=none},
  grid=major, grid style={gray!20}]
\addplot[only marks, mark=*, mark size=0.55pt, blue!70!black] table[x=alpha,y=phi]{phi_beta1p0.dat};
\addlegendentry{$\beta=1$: exact values at $\alpha=r/s$, $s\le20$}
\addplot[only marks, mark=*, mark size=0.55pt, red!70!black] table[x=alpha,y=phi]{phi_beta2p5.dat};
\addlegendentry{$\beta=2.5$: exact values at $\alpha=r/s$, $s\le20$}
\addplot[dashed, blue!60] table[x=alpha,y=low]{lund_beta1p0.dat};
\addlegendentry{Lundberg lower bound $e^{-\theta}$}
\addplot[dashed, red!60, forget plot] table[x=alpha,y=low]{lund_beta2p5.dat};
\addplot[thick, blue!70!black, forget plot] coordinates {(0,1) (1,1)};
\addplot[thick, red!70!black, forget plot] coordinates {(0,1) (2.5,1)};
\addplot[dotted, thick, gray, forget plot] coordinates {(0,0.5) (5,0.5)};
\node[font=\scriptsize, gray] at (axis cs:4.55,0.52) {$\alpha\to\infty$: $b=\tfrac12$};
\draw[gray, ->] (axis cs:3.2,0.70) -- (axis cs:2.53,0.985);
\node[font=\scriptsize, gray!40!black, align=center, fill=white, inner sep=1pt] at (axis cs:3.45,0.66) {threshold $\alpha_\beta=\beta=2.5$\\(earlier bound: $2\le\alpha_\beta\le3$)};
\draw[gray, ->] (axis cs:0.6,0.80) -- (axis cs:0.98,0.985);
\node[font=\scriptsize, gray!40!black, fill=white, inner sep=1pt] at (axis cs:0.6,0.77) {$\alpha_1=1$};
\end{axis}
\end{tikzpicture}
\caption{The crossing probability $\Phi_\beta(\alpha,0)$ as a function of the slope, computed exactly from Theorem~\ref{thm:C} at all rationals $\alpha=r/s\in(\beta,5]$ with $s\le20$. The function equals $1$ exactly on $[0,\beta]$ and has a jump at every rational $\alpha>\beta$; the lower bound $e^{-\theta}$ of Theorem~\ref{thm:D} is attained at integer slopes.}
\label{fig:phi}
\end{figure}

\subsection{Relation to the literature}\label{sec:lit}
We have tried to separate carefully what is classical from what is new.

\emph{Theorem~\ref{thm:A}} is, in other languages, classical. The level process $S_k=Y_k-pX_k$ is a random walk with increments $+1$ and $-p$, so it is \emph{upward skip-free}, and the geometric form $\rho^{d+1}$ together with the first-step equation \eqref{eq:rho} is the textbook solution of the corresponding gambler's ruin problem \cite[Ch.~XIV]{Feller1}. The series for $\rho$ is Lambert's generalized binomial series \cite[\S5.4]{GKP}, the law of the first crossing is an instance of the hitting time theorem \cite{Kemperman,HofstadKeane}, and the identification with the extinction probability of a Galton--Watson tree goes back to the \L ukasiewicz coding of trees \cite{LeGall}. The Fuss--Catalan enumeration \eqref{eq:ballot} of paths below a line of integer slope is due to Fuss \cite{Fuss}; see \cite{HiltonPedersen} and, for slopes $m/n$, Grossman's formula proved by Bizley \cite{Bizley}. What we add is a unified elementary treatment with three short proofs and the mean location $k/(\beta-p)$ of the first crossing.

The enumerative identity for $F(p,n)$ is a Bizley-type ``one contact'' count \cite{Bizley}; our proof via the Green's function and Lagrange--B\"urmann inversion is short and self-contained.

\emph{Theorem~\ref{thm:C}} is an application of the kernel method of Knuth \cite[Ex.~2.2.1-4]{Knuth} and Banderier--Flajolet \cite{BanderierFlajolet} to a \emph{probabilistic} rather than enumerative question: the walk $sY-rX$ has jumps $+s$ and $-r$ and the distribution of its supremum is governed by the large roots of the kernel. Banderier--Flajolet treat the enumeration of such paths in general, and Banderier--Wallner \cite{BanderierWallner} and Duchon \cite{Duchon} study slopes $2/5$ and $2/3$ in detail. The product formula for the generating function of the supremum is the discrete analogue of the Wiener--Hopf factorization \cite[Ch.~XII]{Feller2}, \cite{Spitzer}; the explicit values such as $\Phi_\beta(1/s,0)=s\beta$ and the Green's function formula for $\Psi_\beta(r/s,0)$ appear to be new in this setting.

\emph{Theorem~\ref{thm:D}} sharpens the earlier sufficient condition $\beta\ge\lceil\alpha\rceil$ of \cite{Javaheri07} to the necessary and sufficient condition $\beta\ge\alpha$, and replaces the estimate $\lfloor\beta\rfloor\le\alpha_\beta\le\lfloor\beta\rfloor+1$ for the critical slope by the exact value $\alpha_\beta=\beta$. Part~(ii) is the Lundberg inequality \cite{Asmussen} combined with the observation that overshoots are at most one; the explicit lower bound with constant $e^{-\theta}$ uses the special structure of monotonic walks. Part~(iii), describing the fine structure of $\alpha\mapsto\Phi_\beta(\alpha,0)$, makes precise and proves the continuity statements made without proof in \cite{Javaheri07}.

\begin{figure}[t]
\centering
\begin{tikzpicture}[scale=1.0]
\newcommand{\panel}[3]{%
  \begin{scope}[xshift=#1]
    \draw[gray!30] (0,0) grid (1,2);
    \draw[gray!70] (-0.2,-0.4) -- (1.15,2.3);
    \node[gray, font=\scriptsize, right] at (1.15,2.3) {$y=2x$};
    \draw[very thick, blue!70!black, -latex] #2;
    \foreach \q in {(0,0),(1,0),(0,1),(1,1),(0,2),(1,2)} \fill[gray] \q circle (1.3pt);
    \node[font=\small, align=center] at (0.5,-0.9) {#3};
  \end{scope}}
\panel{0cm}{(0,0) -- (1,0) -- (1,1) -- (1,2)}{$\mathsf{RUU}$: below}
\panel{4cm}{(0,0) -- (0,1) -- (0,2) -- (1,2)}{$\mathsf{UUR}$: above}
\panel{8cm}{(0,0) -- (0,1) -- (1,1) -- (1,2)}{$\mathsf{URU}$: jumps over}
\end{tikzpicture}
\caption{The $F(2,1)=3$ paths from $(0,0)$ to $(1,2)$ that meet $y=2x$ only at the endpoints. The path $\mathsf{URU}$ starts above the line and jumps below it without touching it; such paths were missed in \cite{Javaheri07}.}
\label{fig:jump}
\end{figure}

\subsection*{Organization}
Section~\ref{sec:prelim} reduces everything to a one-dimensional random walk and proves the basic dichotomy. Section~\ref{sec:integer} treats integer slopes (Theorem~\ref{thm:A}), Section~\ref{sec:returns} the return probabilities (Theorem~\ref{thm:B}), Section~\ref{sec:rational} rational slopes (Theorem~\ref{thm:C}) and Section~\ref{sec:real} arbitrary slopes (Theorem~\ref{thm:D}). Asymptotics are collected in Section~\ref{sec:asym}, and Section~\ref{sec:final} lists open problems. All formulas were also checked numerically (Appendix~\ref{app:num}).

\section{Preliminaries: the level process}\label{sec:prelim}

Throughout, $a=1/(\beta+1)$, $b=\beta/(\beta+1)$, and $\W=(X_k,Y_k)$ is the $\beta$-biased monotonic walk from the origin. For a slope $\alpha\ge0$ define the \emph{level process}
\[
S_k=S_k^{(\alpha)}:=Y_k-\alpha X_k,\qquad k\ge0.
\]
This is a random walk on $\R$ started at $0$ with i.i.d.\ increments $\xi$ satisfying $\Prob(\xi=1)=b$ and $\Prob(\xi=-\alpha)=a$, so that
\begin{equation}\label{eq:drift}
\mu:=\E\xi=\frac{\beta-\alpha}{\beta+1},\qquad \sigma^2:=\operatorname{Var}\xi=\frac{\beta(1+\alpha)^2}{(\beta+1)^2}.
\end{equation}
A lattice point $(x,y)$ lies strictly above the line $y=\alpha x+d$ if and only if $y-\alpha x>d$. Hence, writing $M:=\sup_{k\ge0}S_k\in[0,\infty]$,
\[
\Phi_\beta(\alpha,d)=\Prob(M>d),\qquad
\Psi_\beta(\alpha,d)=\Prob\bigl(S_k=d \text{ for some } k\ge1\bigr).
\]
Two structural features of $S$ will be used repeatedly: the upward jumps all have size exactly $1$, so the \emph{overshoot} over any level is at most $1$; and when $\alpha=p$ is an integer, $S$ lives on $\Z$ and is \emph{upward skip-free}: it cannot pass from below $k$ to above $k$ without visiting $k$.

\begin{proposition}[Basic dichotomy]\label{prop:dichotomy}
Let $\alpha,d\ge0$.
\begin{enumerate}[label=\textup{(\alph*)}]
\item If $\beta\ge\alpha$, then $\limsup_k S_k=+\infty$ almost surely, so $\Phi_\beta(\alpha,d)=1$.
\item If $\beta<\alpha$, then $S_k\to-\infty$ almost surely and $\Phi_\beta(\alpha,d)\le\Phi_\beta(\alpha,0)<1$.
\item $\Phi_\beta(\alpha,d)$ is nonincreasing in $\alpha$ and in $d$, and nondecreasing in $\beta$.
\end{enumerate}
\end{proposition}

\begin{proof}
(a) If $\beta>\alpha$ then $\mu>0$ and $S_k/k\to\mu$ by the strong law of large numbers. If $\beta=\alpha$ then $\mu=0$ and $\xi$ is not a.s.\ constant, and a mean-zero nondegenerate random walk satisfies $\limsup_k S_k=+\infty$ a.s.\ \cite[Thm.~4.1.2]{Durrett}.

(b) Now $\mu<0$, so $S_k\to-\infty$ and $M<\infty$ a.s. Suppose $\Prob(M>0)=1$. Let $\tau_1=\inf\{k:S_k>0\}$ and inductively $\tau_{j+1}=\inf\{k>\tau_j:S_k>S_{\tau_j}\}$ (strict ascending ladder epochs). By the strong Markov property each $\tau_j$ is a.s.\ finite and the ladder heights $S_{\tau_{j+1}}-S_{\tau_j}$ are i.i.d.\ and positive, so $S_{\tau_j}\to\infty$, contradicting $M<\infty$. Hence $\Prob(M>0)<1$.

(c) Monotonicity in $\alpha$ and $d$ is clear from the definition. For $\beta$, couple the walks with a single i.i.d.\ uniform sequence $(U_k)$, taking an up-step at time $k$ iff $U_k<\beta/(\beta+1)$. Increasing $\beta$ turns some right-steps into up-steps; this can only increase $Y_k-\alpha X_k$ for every $k$.
\end{proof}

In particular, the critical slope $\alpha_\beta:=\sup\{\alpha:\Phi_\beta(\alpha,0)=1\}$ equals $\beta$, and $\Phi_\beta(\alpha,0)=1$ holds on the closed interval $[0,\beta]$. This proves Theorem~\ref{thm:D}(i).

\section{Integer slopes: crossing probabilities}\label{sec:integer}

\subsection{Fuss--Catalan numbers and the generalized binomial series}
For real $t\ge1$ let
\begin{equation}\label{eq:Bt}
\mathcal B_t(w)=\sum_{n\ge0}\frac{1}{(t-1)n+1}\binom{tn}{n}w^n
=\sum_{n\ge0}\frac{1}{tn+1}\binom{tn+1}{n}w^n,
\end{equation}
Lambert's generalized binomial series \cite[(5.58)]{GKP}; here $\binom{x}{n}=x(x-1)\cdots(x-n+1)/n!$. By Stirling's formula the radius of convergence is
\[
R_t=\frac{(t-1)^{t-1}}{t^t},
\]
and since the coefficients are $O(n^{-3/2}R_t^{-n})$ the series also converges absolutely on $|w|=R_t$. The function $\zeta\mapsto\zeta(1-\zeta)^{t-1}$ increases on $[0,1/t]$ from $0$ to $R_t$.

\begin{proposition}\label{prop:binomial}
For every real $t\ge1$ and every $z\in[0,1/t]$,
\begin{equation}\label{eq:binomial}
\mathcal B_t\bigl(z(1-z)^{t-1}\bigr)=\frac{1}{1-z}.
\end{equation}
\end{proposition}

This is classical \cite[(5.59)]{GKP} and follows from Lagrange inversion. 

\begin{proof}
Write $\alpha=t-1\ge0$ and $c_n=\frac{1}{n\alpha+n+1}\binom{n\alpha+n+1}{n}$; \eqref{eq:binomial} is equivalent to
\begin{equation}\label{eq:binomial2}
L(z):=\sum_{n\ge0}c_n z^n(1-z)^{n\alpha+1}=1 .
\end{equation}
For $\alpha=0$ this is the geometric series. Let $\alpha>0$.

\emph{Step 1: a polynomial identity.} For every $k\ge1$ and every real $\alpha$,
\begin{equation}\label{eq:poly}
\sum_{n=0}^k(-1)^n\binom{k}{n}\binom{n\alpha+n}{k-1}=0 .
\end{equation}
The left side is a polynomial in $\alpha$, so it suffices to verify \eqref{eq:poly} for all even integers $\alpha=m\ge2$. Put $J(y)=\sum_{n=0}^k\binom{k}{n}y^{n(m+1)}=(1+y^{m+1})^k$. Since $m+1$ is odd, $y=-1$ is a root of $1+y^{m+1}$, hence a root of $J$ of multiplicity at least $k$, and $J^{(k-1)}(-1)=0$. Differentiating term by term,
\[
0=J^{(k-1)}(-1)=(k-1)!\sum_{n=0}^k\binom{k}{n}\binom{n(m+1)}{k-1}(-1)^{n(m+1)-k+1},
\]
and $(-1)^{n(m+1)}=(-1)^n$ because $m$ is even. This is \eqref{eq:poly}.

\emph{Step 2: the identity for small $z$.} Expand $(1-z)^{n\alpha+1}=\sum_{l\ge0}(-1)^l\binom{n\alpha+1}{l}z^l$. For real $x\ge0$ we have $|x-j|\le x+j$, hence $\bigl|\binom{x}{l}\bigr|\le\binom{x+l-1}{l}=[z^l](1-z)^{-x}$ and $\sum_l\bigl|\binom{x}{l}\bigr|z^l\le(1-z)^{-x}$ for $0\le z<1$. Therefore, if $z\ge0$ is so small that
\begin{equation}\label{eq:small}
z(1-z)^{-\alpha}<R_{\alpha+1},
\end{equation}
then $\sum_{n,l}c_n\bigl|\binom{n\alpha+1}{l}\bigr|z^{n+l}\le\sum_n c_n z^n(1-z)^{-n\alpha-1}<\infty$, and the double series may be rearranged. A direct computation with Gamma functions gives, for $k\ge1$,
\[
\sum_{n+l=k}(-1)^l c_n\binom{n\alpha+1}{l}=\frac{(-1)^k}{k}\sum_{n=0}^k(-1)^n\binom kn\binom{n\alpha+n}{k-1},
\]
which vanishes by \eqref{eq:poly}. Thus $L(z)=1$ for all $z\ge0$ satisfying \eqref{eq:small}.

\emph{Step 3: analytic continuation.} Fix $c<1/(\alpha+1)$. The map $\omega(z)=z(1-z)^\alpha$ is holomorphic near $[0,c]$ and $|\omega|<R_{\alpha+1}$ on $[0,c]$, hence on a connected complex neighbourhood $U$ of $[0,c]$. On $U$, $L(z)=(1-z)\mathcal B_{\alpha+1}(\omega(z))$ is holomorphic and equals $1$ on the real interval of Step~2. By the identity theorem $L\equiv1$ on $U$, in particular on $[0,c]$. Since $c<1/(\alpha+1)$ was arbitrary, \eqref{eq:binomial2} holds on $[0,1/(\alpha+1))$, and at $z=1/(\alpha+1)$ it follows by Abel's theorem, the series converging on the boundary of its disc.
\end{proof}

\begin{remark}
A shorter route is Lagrange inversion: \eqref{eq:binomial} says that $1/(1-z(w))$, where $z(w)$ is the power series solution of $z=w(1-z)^{-(t-1)}$, has the coefficients \eqref{eq:Bt}; see \eqref{eq:powers} below with $k=1$.
\end{remark}

We shall also need the powers of $\mathcal B_t$ and a companion series \cite[(5.60)--(5.61)]{GKP}: for real $t\ge1$, integer $k\ge1$ and $|w|<R_t$,
\begin{equation}\label{eq:powers}
\mathcal B_t(w)^k=\sum_{n\ge0}\frac{k}{tn+k}\binom{tn+k}{n}w^n,\qquad
\sum_{n\ge0}\binom{tn}{n}w^n=\frac{1}{1-tz},
\end{equation}
where $w=z(1-z)^{t-1}$, $0\le z<1/t$. Both follow from the Lagrange--B\"urmann formula \cite[Thm.~A.2]{FlajoletSedgewick}: if $z=w\varphi(z)$ with $\varphi(z)=(1-z)^{-(t-1)}$, then $[w^n]\,\psi(z)=\frac1n[\zeta^{n-1}]\psi'(\zeta)\varphi(\zeta)^n$ and $\sum_n w^n[\zeta^n]\psi(\zeta)\varphi(\zeta)^n=\psi(z)/(1-w\varphi'(z))$. Taking $\psi(\zeta)=(1-\zeta)^{-k}$ gives the first identity; taking $\psi(\zeta)=(1-\zeta)^{-1}$ and using $w\varphi'(z)=(t-1)z/(1-z)$ gives the second.

For an integer $p\ge1$ set $H_p:=\mathcal B_{p+1}$, so that $H_p(x)=\sum_n M(p,n)x^n$ with $M(p,n)=\frac{1}{pn+1}\binom{(p+1)n}{n}$, and $R_{p+1}=p^p/(p+1)^{p+1}$.

\begin{corollary}\label{cor:H}
Let $p\ge1$ and $0\le x\le R_{p+1}$. Then $H_p(x)$ is the smallest positive root of $xy^{p+1}=y-1$, and $H_p(x)=1/(1-z(x))$ where $z(x)$ is the unique root of $z(1-z)^p=x$ in $[0,1/(p+1)]$.
\end{corollary}

\begin{proof}
By Proposition~\ref{prop:binomial}, $H_p(x)=1/(1-z(x))$, and then $xH_p^{p+1}=z/(1-z)=H_p-1$. Let $y>0$ be any root of $xy^{p+1}=y-1$. Then $y\ge1$, since $x\ge0$. Put $z'=1-1/y\in[0,1)$, so that $z'(1-z')^p=x$. If $z'\le1/(p+1)$ then $z'=z(x)$ by injectivity; otherwise $z'>1/(p+1)\ge z(x)$. In both cases $y=1/(1-z')\ge1/(1-z(x))=H_p(x)$.
\end{proof}

\subsection{Counting paths below a line}\label{sec:counting}
Fix an integer $p\ge1$. For integers $n,d\ge0$ let $S(n,d)$ be the number of monotonic lattice paths from $(0,0)$ to $(n,pn+d)$ that stay weakly below the line $y=px+d$; in level language these are words in $\{\mathsf U,\mathsf R\}$ with $n$ letters $\mathsf R$ and $pn+d$ letters $\mathsf U$ whose level never exceeds $d$. Thus $S(n,0)=M(p,n)$ \cite{HiltonPedersen}.

\begin{lemma}\label{lem:ballot}
For $n,d\ge0$,
\begin{equation}\label{eq:ballot}
S(n,d)=\frac{d+1}{(p+1)n+d+1}\binom{(p+1)n+d+1}{n},\qquad
\sum_{n\ge0}S(n,d)x^n=H_p(x)^{d+1}.
\end{equation}
\end{lemma}

\begin{proof}
We prove $S(n,d+1)=\sum_{i=0}^n S(i,0)\,S(n-i,d)$. A path counted by $S(n,d+1)$ ends at level $d+1\ge1$, so it reaches level $1$; consider the \emph{first} step that brings it to level $1$. That step is an up-step from a point at level $0$, say in column $i$. The part before it goes from level $0$ to level $0$ without exceeding $0$, which is counted by $S(i,0)$; the part after it goes from level $1$ to level $d+1$ without exceeding $d+1$, which after a vertical shift by one unit is counted by $S(n-i,d)$. The decomposition is clearly bijective. 

In terms of generating functions, $\sum_nS(n,d+1)x^n=H_p(x)\sum_nS(n,d)x^n$, and induction on $d$ gives the second formula in \eqref{eq:ballot}. The first follows from \eqref{eq:powers} with $t=p+1$, $k=d+1$.
\end{proof}

The closed formula in \eqref{eq:ballot} is the generalized ballot number; for $d=0$ it gives the Fuss--Catalan numbers, and a bijective proof follows from the cycle lemma of Dvoretzky--Motzkin \cite{DvoretzkyMotzkin}.

\subsection{Proof of Theorem~\ref{thm:A}}
Since $S$ is integer-valued for integer $p$, $\Phi_\beta(p,d)=\Prob(M\ge\lfloor d\rfloor+1)$, so we may assume $d\in\Z_{\ge0}$ and put $k=d+1$. For $p=0$ the walk $S=Y$ is nondecreasing, tends to infinity and $\rho=1$ is the only positive root of \eqref{eq:rho}; so assume $p\ge1$. We give three proofs, each of which illuminates a different aspect.

\begin{proof}[First proof (probabilistic)]
Let $\tau_j=\inf\{k\ge0:S_k=j\}$. Because $S$ is upward skip-free, $\{M\ge k\}=\{\tau_k<\infty\}$, and by the strong Markov property and translation invariance $\Prob(\tau_k<\infty)=\rho^k$ with $\rho:=\Prob(\tau_1<\infty)$. Conditioning on the first step,
\[
\rho=b+a\rho^{p+1},
\]
which is \eqref{eq:rho} divided by $\beta+1$. Let $f(y)=ay^{p+1}-y+b$ and let $y_*$ be any positive root. Put $\rho_N=\Prob(\tau_1\le N)$, so $\rho_N\uparrow\rho$. Reaching level $1$ from level $-p$ within $N$ steps requires $p+1$ successive passages, each of length at most $N$, so $\rho_{N+1}\le b+a\rho_N^{p+1}$. Since $\rho_0=0\le y_*$, induction gives $\rho_N\le b+ay_*^{p+1}=y_*$ for all $N$. Thus $\rho$ is the smallest positive root.

Finally, $f$ is convex on $[0,\infty)$ with $f(0)=b>0$, $f(1)=0$ and $f'(1)=a(p+1)-1=(p-\beta)/(\beta+1)$. If $\beta\ge p$ then $f'(1)\le0$ and $f>0$ on $[0,1)$, so $\rho=1$; if $\beta<p$ then $f(1-\eps)<0$ for small $\eps>0$ and $\rho<1$.
\end{proof}

\begin{proof}[Second proof (enumerative)]
The walk first crosses $y=px+d$ at a point $(n,pn+k)$ precisely when its first $(p+1)n+d$ steps form a path counted by $S(n,d)$ and the next step is an up-step. Hence, by Lemma~\ref{lem:ballot},
\begin{equation}\label{eq:firstcross}
\Prob(\text{first crossing in column } n)=S(n,d)\,a^nb^{pn+k}=\frac{k}{(p+1)n+k}\binom{(p+1)n+k}{n}a^nb^{pn+k},
\end{equation}
and summing over $n$,
\[
\Phi_\beta(p,d)=b^{k}H_p(ab^p)^{k}=\bigl(bH_p(ab^p)\bigr)^{k}.
\]
Note $x:=ab^p=a(1-a)^p\le R_{p+1}$. By Corollary~\ref{cor:H}, $Y=H_p(x)$ is the smallest positive root of $xY^{p+1}=Y-1$, and the substitution $y=bY$ transforms this equation into $ay^{p+1}-y+b=0$, i.e.\ into \eqref{eq:rho}; the substitution preserves order, so $bH_p(ab^p)=\rho$. This gives $\rho=\sum_nM(p,n)a^nb^{pn+1}$. Moreover $z(x)=a$ if $a\le1/(p+1)$, i.e.\ $\beta\ge p$, and then $\rho=b/(1-a)=1$; this is Proposition~\ref{prop:binomial} in disguise.
\end{proof}

\begin{proof}[Third proof (branching)]
Consider a Galton--Watson tree in which each individual has $0$ children with probability $b$ and $p+1$ children with probability $a$. Explore the tree depth first and let $L_j$ be the number of individuals discovered but not yet explored after $j$ steps, minus one; then $(L_j)$ is a random walk with increments $-1$ (prob.\ $b$) and $+p$ (prob.\ $a$) started at $0$, and the tree is finite iff $L$ hits $-1$ (the \L ukasiewicz correspondence \cite[\S1]{LeGall}). Since $-L$ has the law of $S$, the extinction probability of the tree equals $\Prob(\tau_1<\infty)=\rho$, and the probability that the forest grown from $k$ ancestors is finite is $\rho^k$. The classical theory \cite{AthreyaNey} now gives that $\rho$ is the smallest root in $[0,1]$ of $\rho=G(\rho)$, $G(s)=b+as^{p+1}$, and $\rho=1$ iff the mean offspring $(p+1)a\le1$, i.e.\ $\beta\ge p$. Paths counted by $M(p,n)$ correspond to $(p+1)$-ary trees with $n$ internal vertices, recovering one of the interpretations of \cite{HiltonPedersen}.
\end{proof}

\begin{remark}\label{rem:risk}
In the language of risk theory, set $V_k=d-S_k$. Then $V$ decreases by one unit per period and jumps up by $p$ with probability $a$, i.e.\ it is the surplus of a dual (``gains'') risk model with periodic cost $1$ and gains of size $p+1$ \cite{AGS07}; ruin means $V_k<0$, i.e.\ $S_k>d$. Thus $\Phi_\beta(p,d)$ is the ruin probability with initial capital $d$, and Theorem~\ref{thm:A} gives it in closed form.
\end{remark}

\begin{corollary}[Mean location of the first crossing]\label{cor:mean}
If $\beta>p$, the first crossing of $y=px+d$ occurs in a column $X_\tau$ with $\E X_\tau=(\lfloor d\rfloor+1)/(\beta-p)$. If $\beta=p$, then $\E X_\tau=\infty$.
\end{corollary}

\begin{proof}
Let $k=\lfloor d\rfloor+1$ and $\tau=\tau_k$. For $\beta>p$ the drift $\mu=(\beta-p)/(\beta+1)$ is positive, so $\E\tau<\infty$ and Wald's identity gives $k=\E S_\tau=\mu\E\tau$. Since $X_\tau$ is the number of right-steps up to time $\tau$, Wald's identity again gives $\E X_\tau=a\E\tau=ak/\mu=k/(\beta-p)$. If $\beta=p$, then $\E\tau=\infty$ (otherwise Wald's identity would give $k=\E S_\tau=0$); and $\tau=X_\tau+Y_\tau=(p+1)X_\tau+k$, so $\E X_\tau=\infty$.
\end{proof}

\begin{example}
For $p=1$, \eqref{eq:rho} factors as $(y-1)(y-\beta)$, so $\Phi_\beta(1,d)=\min\{1,\beta\}^{\lfloor d\rfloor+1}$. For $p=2$, $(y-1)(y^2+y-\beta)$ gives $\rho=\min\{1,\tfrac{\sqrt{1+4\beta}-1}{2}\}$; for $\beta=1$ this is the inverse golden ratio $0.618\ldots$. In general \eqref{eq:rho} factors as $(y-1)(y^p+\cdots+y-\beta)$, so for $\beta<p$, $\rho$ is the unique positive root of $y+y^2+\cdots+y^p=\beta$.
\end{example}

\section{Integer slopes: returning to the line}\label{sec:returns}

For $n\ge1$ let $N(p,n)$ be the number of monotonic paths from $(0,0)$ to $(n,pn)$ lying strictly below $y=px$ except at their endpoints, and $F(p,n)$ the number of such paths meeting $y=px$ only at their endpoints (on either side). Let $G_p(x)=\sum_{n\ge1}N(p,n)x^n$.

\begin{lemma}\label{lem:N}
$H_p=1+G_pH_p$, so $G_p(x)=1-1/H_p(x)=z(x)$, the root of $z(1-z)^p=x$ in $[0,1/(p+1)]$. Consequently
\[
N(p,n)=\frac1n\binom{(p+1)n-2}{n-1}\qquad(n\ge1).
\]
\end{lemma}

\begin{proof}
For $n\ge1$, a path counted by $M(p,n)$ returns to level $0$ for the first time at some column $m\in\{1,\dots,n\}$; the part before is strictly below the line except at its endpoints. Hence $M(p,n)=\sum_{m=1}^nN(p,m)M(p,n-m)$ for $n\ge1$, i.e.\ $H_p=1+G_pH_p$. Then $G_p=1-1/H_p=z(x)$ by Corollary~\ref{cor:H}. Since $z=x(1-z)^{-p}$, Lagrange inversion gives $[x^n]z=\frac1n[\zeta^{n-1}](1-\zeta)^{-pn}=\frac1n\binom{pn+n-2}{n-1}$.
\end{proof}

\begin{proof}[Proof of Theorem~\ref{thm:B}]
\emph{Case $d\ge1$.} Since $S$ is upward skip-free and starts at $0<d$, it meets level $d$ iff it reaches level $\ge d$; so $\Psi_\beta(p,d)=\Prob(\tau_d<\infty)=\rho^d$.

\emph{Case $d=0$.} Let $V$ be the total number of visits of $S$ to $0$, including time $0$. By the strong Markov property $V$ is geometric: $\Prob(V\ge j+1)=\Psi^j$ with $\Psi=\Psi_\beta(p,0)$, so $\E V=1/(1-\Psi)$ (with $\E V=\infty$ iff $\Psi=1$). On the other hand $S_k=0$ requires $k=(p+1)n$ with $n$ right-steps, so by \eqref{eq:powers} with $t=p+1$ and $x=ab^p$,
\[
\E V=\sum_{n\ge0}\binom{(p+1)n}{n}a^nb^{pn}=\frac{1}{1-(p+1)z_*},\qquad z_*=z(ab^p).
\]
(When $\beta=p$ we have $z_*=1/(p+1)$, and the series diverges, consistent with recurrence.) Hence $\Psi_\beta(p,0)=(p+1)z_*$. By Corollary~\ref{cor:H} and the second proof of Theorem~\ref{thm:A}, $z_*=1-1/H_p(ab^p)=1-b/\rho$. If $\beta\ge p$ then $\rho=1$ and $z_*=a$, so $\Psi_\beta(p,0)=(p+1)/(\beta+1)$, which equals $1$ iff $\beta=p$. If $\beta<p$ then $z_*<1/(p+1)$ and $\Psi_\beta(p,0)<1$.

\emph{The count $F(p,n)$.} Each first-return path with $n$ right-steps has probability $a^nb^{pn}$, so $\Psi_\beta(p,0)=\sum_{n\ge1}F(p,n)x^n$ with $x=ab^p$. We have shown $\Psi_\beta(p,0)=(p+1)z(x)=(p+1)G_p(x)$. Both sides are power series in $x$ that agree for all $x=a(1-a)^p$ with $a\in(0,1)$, i.e.\ for all $x\in(0,R_{p+1}]$; hence $F(p,n)=(p+1)N(p,n)$, and Lemma~\ref{lem:N} finishes the proof.
\end{proof}

\begin{remark}[Why the factor is $p+1$]
Rotation by $180^\circ$ about the midpoint of the segment from $(0,0)$ to $(n,pn)$ shows that exactly $N(p,n)$ first-return paths lie strictly above the line. The remaining $(p-1)N(p,n)$ first-return paths start above the line and jump below it, which is possible only for $p\ge2$ (Figure~\ref{fig:jump}). Table~\ref{tab:F} lists the first values. A direct bijective explanation of the factor $p+1$, e.g.\ via cyclic shifts of the $(p+1)n$ letters as in the cycle lemma, would be of interest.
\end{remark}

\begin{table}[t]
\centering
\begin{tabular}{c|rrrrr}
\toprule
$p\backslash n$ & 1 & 2 & 3 & 4 & 5\\
\midrule
1 & 2 & 2 & 4 & 10 & 28\\
2 & 3  & 6  & 21  & 90 & 429 \\
3 & 4  & 12  & 60  & 364  & 2448 \\
\bottomrule
\end{tabular}
\vspace{0.5em}
\caption{The number $F(p,n)$ of paths from $(0,0)$ to $(n,pn)$ meeting $y=px$ only at the endpoints, confirmed by exhaustive enumeration.}
\label{tab:F}
\end{table}

\begin{corollary}[Number of visits]
For $\beta\ne p$ the number of returns of $\W$ to the line $y=px$ is geometric with mean $\Psi/(1-\Psi)=(p+1)z_*/(1-(p+1)z_*)$; for $\beta>p$ this equals $(p+1)/(\beta-p)$. For $\beta=p$ the walk returns infinitely often almost surely.
\end{corollary}

\section{Rational slopes}\label{sec:rational}

Let $\alpha=r/s$ with $\gcd(r,s)=1$. The scaled level process $\Sigma_k:=sS_k=sY_k-rX_k$ is a random walk on $\Z$ with jumps $+s$ (prob.\ $b$) and $-r$ (prob.\ $a$). For $s\ge2$ it is \emph{not} upward skip-free, so the crossing probability is no longer a pure power of a single number. Since $\gcd(r,s)=1$, the walk $\Sigma$ is irreducible on $\Z$. A lattice point lies strictly above $y=\frac rs x+d$ iff $\Sigma>sd$, iff $\Sigma\ge\lfloor sd\rfloor+1$. Put
\[
K(z)=bz^{r+s}-z^r+a,
\]
the \emph{kernel} of the walk: $K(z)=z^r\bigl(bz^s+az^{-r}-1\bigr)$.

\begin{lemma}\label{lem:roots}
Assume $\beta<r/s$. Then $K$ has exactly $r$ roots (with multiplicity) in the closed unit disc, among them the simple root $z=1$ and no other root on $|z|=1$, and exactly $s$ roots $v_1,\dots,v_s$ in $|z|>1$. Among the latter there is a simple real root $v_*>1$, and every other root satisfies $|v_j|>v_*$.
\end{lemma}

\begin{proof}
Let $g(t)=t^r-bt^{r+s}-a$ for $t>0$. Then $g(1)=0$ and $g'(1)=r-b(r+s)=(\beta+1)^{-1}(r-\beta s)>0$, while $g(t)\to-\infty$; by Descartes' rule of signs $g$ has at most two positive roots, so it has exactly one other positive root $v_*>1$, and $g>0$ on $(1,v_*)$. For $|z|=t\in(1,v_*)$ we get $|bz^{r+s}+a|\le bt^{r+s}+a<t^r=|z^r|$, so by Rouch\'e's theorem $K$ and $-z^r$ have the same number $r$ of zeros in $|z|<t$; letting $t\downarrow1$ shows that $K$ has $r$ zeros in $|z|\le1$ and none in $1<|z|<v_*$. If $|z|=1$ and $K(z)=0$, then $1=|z^r|=|bz^{r+s}+a|\le1$, forcing $z^{r+s}=1$ and $z^r=1$, hence $z^s=1$, and $z=1$ since $\gcd(r,s)=1$. Similarly, if $|z|=v_*$ and $K(z)=0$, equality in $|bz^{r+s}+a|\le bv_*^{r+s}+a=v_*^r$ forces $z^{r+s}>0$ and then $z^r>0$, so $z=v_*$. The roots $1$ and $v_*$ are simple because $g'(1)\ne0\ne g'(v_*)$.
\end{proof}

\begin{theorem}\label{thm:rational}
Let $\gcd(r,s)=1$, $\beta<r/s$, and let $v_1,\dots,v_s$ be the roots of $K$ in $|z|>1$. The supremum $M_\Sigma=\sup_k\Sigma_k$ satisfies
\begin{equation}\label{eq:WH}
\E\bigl[z^{M_\Sigma}\bigr]=\prod_{j=1}^s\frac{v_j-1}{v_j-z},\qquad|z|\le1 .
\end{equation}
Consequently, with $k=\lfloor sd\rfloor+1$,
\[
\Phi_\beta(r/s,d)=\Prob(M_\Sigma\ge k)=1-[z^{k-1}]\,\frac{1}{1-z}\prod_{j=1}^s\frac{v_j-1}{v_j-z},
\]
and if the $v_j$ are distinct,
\begin{equation}\label{eq:partialfrac}
\Phi_\beta(r/s,d)=\sum_{j=1}^s\Bigl(\prod_{i\ne j}\frac{v_i-1}{v_i-v_j}\Bigr)v_j^{-k}.
\end{equation}
In particular $\Phi_\beta(r/s,0)=1-\prod_{j=1}^s(1-1/v_j)$, and $\Phi_\beta(r/s,d)\sim c_*\,v_*^{-\lfloor sd\rfloor-1}$ as $d\to\infty$, with $c_*>0$.
\end{theorem}

\begin{proof}
Let $\pi_m=\Prob(M_\Sigma=m)$ and $F(z)=\E z^{M_\Sigma}$, which is continuous on $|z|\le1$ and holomorphic inside. Since $(\Sigma_{k+1}-\Sigma_1)_k$ is independent of the first jump $\xi$ and has the law of $\Sigma$, Lindley's equation $M_\Sigma\overset{d}{=}\max(0,\xi+M_\Sigma')$ holds with $M'_\Sigma$ an independent copy. Therefore
\[
F(z)=F(z)\bigl(bz^s+az^{-r}\bigr)+a\sum_{m=0}^{r-1}\pi_m\bigl(1-z^{m-r}\bigr),
\]
because $\xi+M'_\Sigma<0$ only when $\xi=-r$ and $M'_\Sigma<r$. Multiplying by $-z^r$,
\[
F(z)K(z)=-a\,Q(z),\qquad Q(z):=\sum_{m=0}^{r-1}\pi_m\bigl(z^r-z^m\bigr),
\]
a polynomial of degree $\le r$, not identically zero because $\pi_0=1-\Phi_\beta(r/s,0)>0$ by Proposition~\ref{prop:dichotomy}(b). Every zero of $K$ in $|z|\le1$ is a zero of $Q$ of at least the same multiplicity (on the boundary only $z=1$ occurs, which is simple). By Lemma~\ref{lem:roots} there are $r$ such zeros $u_1,\dots,u_r$, so $Q=c\prod_i(z-u_i)$, and since $K=b\prod_i(z-u_i)\prod_j(z-v_j)$, we get $F(z)=C\prod_j(z-v_j)^{-1}$. The normalization $F(1)=1$ gives \eqref{eq:WH}. Writing $\Prob(M_\Sigma\ge k)=1-\sum_{m<k}\pi_m$ gives the coefficient formula. For \eqref{eq:partialfrac}, expand $F(z)=\sum_jA_j/(v_j-z)$ with $A_j=\prod_i(v_i-1)/\prod_{i\ne j}(v_i-v_j)$; then $\pi_m=\sum_jA_jv_j^{-m-1}$ and $\Prob(M_\Sigma\ge k)=\sum_jA_jv_j^{-k}/(v_j-1)$. The asymptotics follow from Lemma~\ref{lem:roots}, since $v_*$ is the unique root of smallest modulus, and $c_*>0$ by the lower bound of Theorem~\ref{thm:D}(ii).
\end{proof}

For $s=1$ the formula reduces to $\Phi_\beta(p,d)=v_*^{-(d+1)}$, and indeed $\rho=1/v_*$ satisfies \eqref{eq:rho}. When $\beta\ge r/s$, Proposition~\ref{prop:dichotomy} gives $\Phi_\beta(r/s,d)=1$.

\begin{example}[Slopes $1/s$]\label{ex:1s}
Let $r=1$. Then $K(z)=(z-1)P(z)$ with $P(z)=b(z^s+\cdots+z)-a$, the $v_j$ are the roots of $P$, and $\prod_j(v_j-1)=(-1)^sP(1)/b$, $\prod_jv_j=(-1)^sP(0)/b$. Hence
\[
\prod_{j=1}^s\Bigl(1-\frac{1}{v_j}\Bigr)=\frac{P(1)}{P(0)}=\frac{bs-a}{-a}=1-s\beta,
\qquad\text{so}\qquad
\Phi_\beta(1/s,0)=\min\{1,s\beta\}.
\]
Thus the probability that $\W$ ever rises strictly above $y=x/s$ is exactly $s\beta$ when $\beta\le1/s$: a strikingly simple answer, given that the walk must first move right and then rise faster than the line. Equivalently, after exchanging the coordinates: a walk with up-bias $\beta'=1/\beta\ge s$ stays weakly above the line $y=sx$ forever with probability $1-s/\beta'$.
\end{example}

\begin{example}[Slope $3/2$, $\beta=1$]
Here $K(z)=\frac12(z-1)(z^4+z^3-z^2-z-1)$, whose quartic factor has two roots outside the unit disc, $v_*\approx1.178724$ and $v_2\approx-1.512876$. Theorem~\ref{thm:rational} yields $\Phi_1(3/2,0)\approx0.748152$, $\Phi_1(3/2,\tfrac12)\approx0.700960$ and $\Phi_1(3/2,2)\approx0.401920$; the Lundberg bounds of Theorem~\ref{thm:D} with $\theta=2\log v_*$ give $0.7197\le\Phi_1(3/2,0)<1$.
\end{example}

\subsection*{Returns for rational slopes}
The line $y=\frac rs x+d$ contains lattice points iff $D:=sd\in\Z$. Let
\[
\mathcal G(D):=\sum_{\substack{x,y\ge0\\ sy-rx=D}}\binom{x+y}{x}a^xb^y
\]
be the expected number of visits of $\Sigma$ to $D$ (the Green's function). The points with $sy-rx=0$ are $(sk,rk)$, so $\mathcal G(0)=\sum_k\binom{(r+s)k}{sk}(b^ra^s)^k$, and $b^ra^s=\beta^r/(\beta+1)^{r+s}$.

\begin{theorem}\label{thm:rationalreturn}
Let $\gcd(r,s)=1$. If $\beta\ne r/s$, then $\mathcal G(0)<\infty$,
\[
\Psi_\beta(r/s,0)=1-\frac{1}{\mathcal G(0)},\qquad
\Psi_\beta(r/s,D/s)=\frac{\mathcal G(D)}{\mathcal G(0)}\quad(D\in\Z_{>0}).
\]
If $\beta=r/s$, then $\Psi_\beta(r/s,D/s)=1$ for every $D\in\Z_{\ge0}$.
\end{theorem}

\begin{proof}
By Stirling's formula the $k$-th term of $\mathcal G(0)$ is asymptotic to $c\,k^{-1/2}\kappa^k$ with $\kappa=\frac{(r+s)^{r+s}}{r^rs^s}b^ra^s$. By the weighted AM--GM inequality, $\kappa\le1$ with equality iff $b=r/(r+s)$, i.e.\ $\beta=r/s$. So $\mathcal G(0)<\infty$ iff $\beta\ne r/s$, which is the transience criterion for the irreducible walk $\Sigma$. For a transient walk, the number of visits to $0$ is geometric with mean $\mathcal G(0)=1/(1-\Psi_\beta(r/s,0))$, and by the strong Markov property and translation invariance $\mathcal G(D)=\Prob(\Sigma\text{ hits }D)\,\mathcal G(0)$. If $\beta=r/s$, then $\Sigma$ is a recurrent irreducible walk on $\Z$, so it visits every integer a.s.
\end{proof}

For $s=1$, \eqref{eq:powers} gives $\mathcal G(0)=1/(1-(p+1)z_*)$, and Theorem~\ref{thm:rationalreturn} reduces to Theorem~\ref{thm:B}. For example $\Psi_1(3/2,0)\approx0.596679$ and $\Psi_2(5/2,0)\approx0.720249$.

\section{Arbitrary real slopes}\label{sec:real}

\begin{proof}[Proof of Theorem~\ref{thm:D}\textup{(ii)}]
Let $\beta<\alpha$ and $\Lambda(\theta)=be^\theta+ae^{-\alpha\theta}$, the moment generating function of $\xi$. It is strictly convex, $\Lambda(0)=1$, $\Lambda'(0)=\mu<0$ and $\Lambda(\theta)\to\infty$ as $\theta\to\infty$, so there is a unique $\theta>0$ with $\Lambda(\theta)=1$. Then $e^{\theta S_k}$ is a martingale. Let $\tau=\inf\{k:S_k>d\}$. For $k<\tau$ we have $S_k\le d$, and since the upward jumps have size $1$, $S_\tau\in(d,d+1]$ on $\{\tau<\infty\}$. Hence $e^{\theta S_{\tau\wedge n}}\le e^{\theta(d+1)}$, and on $\{\tau=\infty\}$, $e^{\theta S_n}\to0$ because $S_n\to-\infty$. Optional stopping and dominated convergence give
\[
1=\E\bigl[e^{\theta S_\tau};\tau<\infty\bigr],
\]
and since $e^{\theta d}<e^{\theta S_\tau}\le e^{\theta(d+1)}$ on $\{\tau<\infty\}$,
\[
e^{-\theta(d+1)}\le\Prob(\tau<\infty)=\Phi_\beta(\alpha,d)<e^{-\theta d}.
\]
If $\alpha=p$ and $d$ are integers, then $S_\tau=d+1$ exactly and equality holds on the left. This gives a fourth proof of Theorem~\ref{thm:A}, with $\rho=e^{-\theta}$ (indeed $be^\theta+ae^{-p\theta}=1$ is \eqref{eq:rho} for $y=e^{-\theta}$).
\end{proof}

\begin{remark}
Since $\Phi_\beta(\alpha,d)$ is nonincreasing in $d$, the bound gives $\Phi_\beta(\alpha,d)=e^{-\theta d+O(1)}$ uniformly in $d$, with the $O(1)$ term in $[-\theta,0]$. For irrational $\alpha$ the increments are non-lattice, and the Cram\'er--Lundberg renewal argument \cite[Ch.~XIII]{Asmussen} sharpens this to $\Phi_\beta(\alpha,d)\sim C e^{-\theta d}$ with $C=\lim_{d\to\infty}\bigl(\E[e^{\theta(S_\tau-d)}\mid\tau<\infty]\bigr)^{-1}\in[e^{-\theta},1]$; for rational $\alpha$ the exact asymptotics are those of Theorem~\ref{thm:rational}. For rational $\alpha=r/s$, $e^{\theta/s}=v_*$.
\end{remark}

\begin{proof}[Proof of Theorem~\ref{thm:D}\textup{(iii)}]
Define $\mathcal R=\sup_{k\ge1}Y_k/X_k\in[0,\infty]$ with the convention $y/0=+\infty$ for $y>0$. The walk visits a point strictly above $y=\alpha x$ iff some ratio $Y_k/X_k$ exceeds $\alpha$, i.e.\ iff $\mathcal R>\alpha$. Thus $\Phi_\beta(\alpha,0)=\Prob(\mathcal R>\alpha)$ is one minus a distribution function, hence right-continuous and nonincreasing, with a jump of size $\Prob(\mathcal R=\alpha)$ at $\alpha$. The event $\{\mathcal R=\infty\}$ is the event that the first step is an up-step, so $\Phi_\beta(\alpha,0)\to b$ as $\alpha\to\infty$.

By the strong law, $Y_k/X_k\to\beta$ a.s., so $\mathcal R\ge\beta$, and by Proposition~\ref{prop:dichotomy}(a), $\Prob(\mathcal R>\beta)=\Phi_\beta(\beta,0)=1$. On the event $\{\beta<\mathcal R<\infty\}$, the supremum of a sequence converging to $\beta$ and exceeding $\beta$ is attained, so $\mathcal R=Y_k/X_k\in\Q$ for some $k$. Hence $\mathcal R$ is a.s.\ rational or infinite, $\Prob(\mathcal R=\alpha)=0$ for irrational $\alpha$, and $\Phi_\beta(\cdot,0)$ is a pure jump function on $(\beta,\infty)$ that is continuous at every irrational point.

Let $\alpha=r/s>\beta$ in lowest terms. With positive probability the walk reaches $(s,r)$ along the path $\mathsf{R}^s\mathsf U^r$, never rising above $y=\alpha x$; from there, by the Markov property, it stays weakly below the line forever with probability $1-\Phi_\beta(\alpha,0)>0$. On this event $\mathcal R=\alpha$, so $\Prob(\mathcal R=\alpha)>0$. Since the rationals are dense, $\Phi_\beta(\cdot,0)$ is strictly decreasing on $[\beta,\infty)$.
\end{proof}

\section{Asymptotics}\label{sec:asym}

\begin{proposition}\label{prop:asym}
\begin{enumerate}[label=\textup{(\alph*)}]
\item \emph{Steep lines.} For fixed $\beta>0$, as $p\to\infty$,
\[
\Phi_\beta(p,0)=b+ab^{p+1}+(p+1)a^2b^{2p+1}+O\bigl(p^3b^{3p}\bigr).
\]
\item \emph{Near-critical slopes.} As $\alpha\downarrow\beta$,
\[
\theta=\frac{2(\alpha-\beta)}{\beta+\alpha^2}+O\bigl((\alpha-\beta)^2\bigr)=\frac{2(\alpha-\beta)}{\beta(\beta+1)}+O\bigl((\alpha-\beta)^2\bigr),
\]
so that $\Phi_\beta(\alpha,d)=\exp\bigl(-\tfrac{2(\alpha-\beta)}{\beta(\beta+1)}d\,(1+o(1))\bigr)$ whenever $d\to\infty$ and $(\alpha-\beta)\to0$; in particular, for integer slopes, $1-\Phi_\beta(p,0)\sim\frac{2(p-\beta)}{p(p+1)}$ as $\beta\uparrow p$.
\end{enumerate}
\end{proposition}

\begin{proof}
(a) By Theorem~\ref{thm:A}, $\Phi_\beta(p,0)=b\sum_nM(p,n)x^n$ with $x=ab^p$, and $M(p,1)=1$, $M(p,2)=p+1$. Since $M(p,n)\le\binom{(p+1)n}{n}\le(e(p+1))^n$, the tail $\sum_{n\ge3}$ is bounded by a geometric series with first term $(e(p+1)ab^p)^3$ and ratio tending to $0$.

(b) Expanding $\Lambda(\theta)=1+\mu\theta+\frac12(b+a\alpha^2)\theta^2+O(\theta^3)$ and solving $\Lambda(\theta)=1$ gives $\theta=-2\mu/(b+a\alpha^2)+O(\mu^2)=2(\alpha-\beta)/(\beta+\alpha^2)+O((\alpha-\beta)^2)$, using $(\beta+1)(b+a\alpha^2)=\beta+\alpha^2$; replacing $\alpha^2$ by $\beta^2$ changes this only at second order. The statement about $\Phi_\beta$ follows from Theorem~\ref{thm:D}(ii). For integer slopes, $\rho=e^{-\theta}$.
\end{proof}

Part (a) replaces the asymptotic statement of \cite[\S4]{Javaheri07} by an explicit expansion whose leading terms have a transparent meaning: $b$ is the probability of crossing with the very first step, $ab^{p+1}$ that of the path $\mathsf{RU}^{p+1}$, and so on. Part (b) is the diffusion (heavy-traffic) approximation: near criticality, $S$ is close to a Brownian motion with drift $\mu$ and variance $\sigma^2$, and $2|\mu|/\sigma^2\to2(\alpha-\beta)/(\beta(\beta+1))$.

\section{Concluding remarks and open problems}\label{sec:final}

\begin{enumerate}[leftmargin=2em]
\item \emph{General monotone walks.} Nothing in Sections~\ref{sec:prelim}, \ref{sec:rational} and~\ref{sec:real} uses the planar structure beyond the reduction to the level process. For a monotone walk in $\Z^m$ with steps $e_1,\dots,e_m$ chosen with probabilities $q_1,\dots,q_m$, crossing the hyperplane $x_m=\sum_{i<m}\alpha_ix_i+d$ is governed by the walk with increments $+1$ and $-\alpha_i$; the threshold is $q_m\ge\sum_{i<m}\alpha_iq_i$, the Lundberg sandwich of Theorem~\ref{thm:D}(ii) holds verbatim, and for rational $\alpha_i$ the kernel method of Section~\ref{sec:rational} applies with the obvious changes.
\item \emph{A bijection.} Find a bijective proof of $F(p,n)=(p+1)N(p,n)$.
\item \emph{Multiple crossings.} For $p\ge2$, the walk may cross $y=px$ from above to below without touching it. Determine the law of the number of upcrossings and downcrossings of the line; the return count of Theorem~\ref{thm:B} is only the first step.
\item \emph{Fine structure in $\alpha$.} The jump of $\Phi_\beta(\cdot,0)$ at a rational $\alpha=r/s>\beta$ equals $\Prob(\mathcal R=r/s)$. How does it depend on $s$ (e.g.\ is it of order $\kappa^s$ for some $\kappa<1$), and what is the Hausdorff-type regularity of the resulting devil's-staircase-like function?
\item \emph{Curved boundaries and varying bias.} Replace the line by $y=\alpha x+c\sqrt{x}$ or $y=\alpha x+c\sqrt{x\log\log x}$ at the critical slope $\alpha=\beta$ (Kolmogorov--Petrovsky-type tests), or let the bias depend on the position, as in chemotaxis along a gradient \cite{Codling}.
\end{enumerate}

\appendix
\section{Numerical verification}\label{app:num}
All closed forms in this paper were checked against independent computations: crossing and hitting probabilities were obtained by solving the linear equations $h(i)=b\,h(i+s)+a\,h(i-r)$ for the level process on a large truncated window ($|i|\le 6000$), and lattice-path counts by exhaustive enumeration. Representative results:

\begin{center}
\small
\begin{tabular}{lcc}
\toprule
quantity & linear system & formula \\
\midrule
$\Phi_{1.3}(2,3)$ & $0.3080361447$ & $\rho^4=0.3080361447$ \\
$\Psi_{1}(2,0)$ & $0.5729490169$ & $3(1-b/\rho)=0.5729490169$ \\
$\Psi_{3}(2,0)$ & $0.7500000000$ & $3/4$ \\
$\Psi_{5}(3,0)$ & $0.6666666667$ & $4/6$ \\
$\Psi_{2}(3,0)$ & $0.7099948573$ & $4(1-b/\rho)=0.7099948573$ \\
$\Phi_{1}(3/2,2)$ & $0.4019198590$ & Thm.~\ref{thm:rational}: $0.4019198590$ \\
$\Phi_{0.05}(1/7,0)$ & $0.3500000000$ & $7\beta=0.35$ \\
$\Psi_{0.5}(2/3,0)$ & $0.7028440076$ & $1-1/\mathcal G(0)=0.7028440076$ \\
\bottomrule
\end{tabular}
\end{center}
\noindent For irrational slopes, $\Phi_1(\sqrt2,2.5)\approx0.4324$ lies between the Lundberg bounds $0.3624$ and $0.4843$, and $\Phi_2(\pi,1)\approx0.6495$ between $0.6348$ and $0.7967$.

\end{document}